\documentclass[twoside,12pt]{article}
\usepackage{CJK}
\usepackage{indentfirst}
\usepackage{bm}
\usepackage[]{hyperref}
\usepackage[numbers,sort&compress]{natbib}
\usepackage{tikz}
\usepackage{float}
\usepackage{color}
\usepackage{graphicx}
\usepackage{subfig}
\usepackage{amsmath}
\usepackage{mathrsfs}
\usepackage{amssymb}
\usepackage{amsthm}
\newcommand{\RNum}[1]{\uppercase\expandafter{\romannumeral #1\relax}}

\makeatletter
\@addtoreset{equation}{section}
\makeatother
\begin{document}
\begin{CJK*}{GBK}{song}
\allowdisplaybreaks


\begin{center}
\LARGE\bf Asymptotic behavior for a new higher-order nonlinear Schr\"{o}dinger equation 
\end{center}
\footnotetext{\hspace*{-.45cm}\footnotesize $^*$Corresponding author: Y.F. Zhang. E-mail: mathzhang@126.com, \\}
\begin{center}
\ \\Hongyi Zhang$^{1}$, Yufeng Zhang$^{1*}$, Binlu Feng$^{2}$
\end{center}
\begin{center}
\begin{small} \sl
{$^{1}$School of Mathematics, China University of Mining and Technology, Xuzhou, Jiangsu, 221116, People's Republic of China.\\
$^{2}$School of Mathematics and Information Sciences, Weifang University, Weifang,Shandong, 261061, People's Republic of China.}
\end{small}
\end{center}
\vspace*{2mm}
\begin{center}
\begin{minipage}{15.5cm}
\parindent 20pt\footnotesize


 \noindent {\bfseries
Abstract} We investigate the Cauchy problem of a new higher-order nonlinear Schr\"{o}dinger equation (NHNSE) with weighted Sobolev initial data which is derived by ourselves. By applying $\bar{\partial}$-steepest descent method, we derive the long-time asymptotics of the NHNSE. Explicit steps are as follows: first of all, based on the spectral analysis of a Lax pair and scattering matrice, the solution of the NHNSE is exhibted through solving the corresponding Riemann-Hilbert problem. Secondly, by applying some properties of the Riemann-Hilbert problem, we obtain the long-time asymptotics of the solution to the NHNSE. As we know that the properties of the NHNSE presented in the paper have not been found in any scholar journals.
\end{minipage}
\end{center}
\begin{center}
\begin{minipage}{15.5cm}
\begin{minipage}[t]{2.3cm}{\bf Keywords}\end{minipage}
\begin{minipage}[t]{13.1cm}
  Long-time asymptotics; Riemann-Hilbert problem; $\bar{\partial}$-steepest descent method.
\end{minipage}\par\vglue8pt
\end{minipage}
\end{center}

\tableofcontents

\section{Introduction}  

Researing solitons and their interactions not only helps to elucidate the motion laws of matter under nonlinear interactions, but also promotes the development of methods and techniques for solving nonlinear evolution equations.
 As a result, a large number of methods have been proposed, developed and promoted \cite{Ablowitz83,Ablowitz91,Deift93,Hirota71,
Hirota04,Konopelchenko93,Matveev91,Rogers82,Rogers02,Sklyanin79,Zakharov74,Zakharov85,Xia15,Zhang22,Chai21}.
One of the most important methods is the well-known inverse scattering transform (IST) method proposed by Gardner, Green, Kruskal and Miura 1967 \cite{Gardner67}, and since then the IST has widely been used in the study of integrable systems and soliton solutions with analyzing long-time asymptotic behavior.

In 1973, Manakov first used the IST method to discuss the long-time asymptotic behavior of the nonlinear wave equations \cite{Manakov73}.
Since then, the IST methods have received widespread attention and development in nonlinear sciences \cite{Zakharov76,Bikbaev88,Fokas92}.
Particularly, in 1993 the remarkable nonlinear steepest descent method was proposed by Deift and Zhou to investigate the long-time asymptotic behavior of integrable equations via the Riemann-Hilbert problem \cite{Deift93}. They successfully applied their nonlinear steepest descent approach to the nonlinear Schr\"{o}dinger equation \cite{Deift03}.
Later in 2016, the Deift-Zhou method was generalized to the so-called $\bar{\partial}$-steepest descent approach \cite{Cuccagna16,Jenkins18,Borghese18,Dieng19,Yang22,Cheng22}.
Those important models for the study of asymptotic stability of N-soliton solutions by using the $\bar{\partial}$-steepest descent approach include: 1) The defocusing nonlinear schr\"{o}dinger equation by Cuccagna and Jenkins \cite{Cuccagna16};
2) The derivative nonlinear schr\"{o}dinger equation by Jenkins, Liu, Perry and Sulem
\cite{Jenkins18};
3) The focusing nonlinear schr\"{o}dinger equation by Borghese, Jenkins, McLaughlin, and Miller \cite{Borghese18}, and other integrable PDE models by Dieng and McLaughlin \cite{Dieng19}
4)The cubic Camassa-Holm equation in space-time solitonic regions by Yang and Fan \cite{Yang22};
and 5)The focusing Fokas-Lenells equation in the solitonic region of space-time by Cheng and Fan \cite{Cheng22}.

Introducing the following Lax pair:
\begin{align}
&\psi_x=-iz\sigma_3\psi+\tilde{P}\psi,\label{lax1a}\\
&\psi_t=[(\alpha z^3+\beta z^2+\gamma z+\delta)\sigma_3+\tilde{Q}]\psi,\label{lax1b}
\end{align}
where\\
$~~~~~~~~~~~~~~~~~~~~~~~~~~~~~~~~~~~~~~~~\tilde{P}=\left(\begin{array}{cc}
0 & q \\
r & 0
\end{array}\right)$\\
and\\
$~~~~~~~~~~~~\tilde{Q}=i \alpha z^{2} \tilde{P}-iz\left(\begin{array}{cc}
\frac{i \alpha}{2} q r & -\frac{i \alpha}{2} q_{x}-\beta q \\
\frac{i \alpha}{2} r_{x}-\beta r & -\frac{i \alpha}{2} q r
\end{array}\right)$\\
$~~~~~~~~~~~~~~~~~-\left(\begin{array}{cc}
\frac{i \alpha}{4}\left(q r_{x}-r q_{x}\right)-\frac{\beta}{2} q r & -\frac{i \alpha}{4}\left(-q_{x x}+2 q^{2} r\right)+\frac{\beta}{2} q_{x}-i \gamma q \\
-\frac{i \alpha}{4}\left(-r_{x x}+2 q r^{2}\right)-\frac{\beta}{2} r_{x}-i \gamma r & -\frac{i \alpha}{4}\left(q r_{x}-r q_{x}\right)+\frac{\beta}{2} q r
\end{array}\right)$,\\
 $\alpha$, $\beta$, $\gamma$ and $\delta$ are pure imaginary numbers.
It is easy to verify that the compatibility conditions of \eqref{lax1a} and \eqref{lax1b} gives rise to
\begin{equation}\label{akns}
\begin{cases}
{q_t}=-\frac{i\alpha}{4}({q_{xxx}}-6qr{q_x})-\frac{\beta}{2}({q_{xx}}-2q^{2}r)+i\gamma {q_x}+2\delta q,\\
{r_t}=-\frac{i\alpha}{4}({r_{xxx}}-6qr{r_x})+\frac{\beta}{2}({r_{xx}}-2qr^{2})+i\gamma {r_x}-2\delta r,\\
\end{cases}
\end{equation}
The advantages of Eq.\eqref{akns} can reduce to some well-known equations which have important physical applications as follows:\\
(I)Taking $\alpha=-4i$, $\beta=\gamma=\delta=0$, $r=-1$, system \eqref{akns} reduces to KdV equation:
\begin{equation}
{q_t}+6q{q_x}+{q_{xxx}}=0.
\end{equation}
(II)Choosing $\alpha=-4i$, $\beta=\gamma=\delta=0$, $r=-q$, system \eqref{akns} reduces to modified KdV equation:
\begin{equation}
{q_t}+6q^{2}{q_x}+{q_{xxx}}=0.
\end{equation}
(III)Let $\beta=-2i$, $\alpha=\gamma=\delta=0$, $r=\mp q^{\ast}$, system \eqref{akns} reduces to nonlinear Schr\"{o}dinger equation
\begin{equation}
i{q_t}+{q_{xx}}\pm2q^{2}q^{\ast}=0.
\end{equation}
Taking $\alpha=i$, $\beta=-i$, $\gamma=i$, $\delta=i$ and $r=q^{*}$, the system \eqref{akns} reduces to
\begin{equation}\label{hnls0a}
 iq_{t}-\frac{i}{4}q_{xxx}+\frac{3i}{2}|q|^{2}q_x+\frac{1}{2}q_{xx}+i q_{x}+2q-|q|^2q=0.
\end{equation}
The Lax pair of \eqref{hnls0a} obviously reads:
\begin{align}
&\psi_x=-iz\sigma_3\psi+P\psi,\label{lax0a}\\
&\psi_t=[(iz^3-iz^2+iz+i)\sigma_3+Q]\psi,\label{lax0b}
\end{align}
where $\psi(z;x,t)$ is a $2\times2$ matrix and\\
$~~~~~~~~~~~~~~~~~~~~~\sigma_3=\left(\begin{array}{cc}
1 & 0 \\
0 & -1
\end{array}\right)$,~~
$P=\left(\begin{array}{cc}
0 & q \\
q^{*} & 0
\end{array}\right)$,\\

$Q=-z^2P+\frac{iz}{2}(P^2+P_x)\sigma_3+\frac{1}{4}[P,P_x]-\frac{1}{2}P^3+\frac{1}{4}P_{xx}
+zP-\frac{i}{2}(P^2+P_{x})\sigma_3-P$.\\
In this paper, we study the Cauchy problem of the NHNSE by us:
\begin{align}
&iq_{t}-\frac{i}{4}q_{xxx}+\frac{3i}{2}|q|^{2}q_x+\frac{1}{2}q_{xx}+i q_{x}+2q-|q|^2q=0,\label{hnls1a}\\
&q(x,0)=q_0(x) \in H^{1,1}(\mathbb{R}), \label{hnls1b}
\end{align}
where $H^{1,1}(\mathbb{R})$ is weighted Sobolev space\\
$~~~~~~~~~~~~~~~~~~~~~~~~~~~H^{1,1}(\mathbb{R})=\{f(x)\in L^2(\mathbb{R}):f^{'}(x),xf(x)\in L^2(\mathbb{R})\}.$\\

\section{Properties of characteristic function and Riemann-Hilbert problem}  
Under the initial condition \eqref{hnls1b}, Lax pair \eqref{lax0a}-\eqref{lax0b} has the following asymptotic form of Jost solution
\begin{equation}
\psi \sim e^{i[-zx+(z^3-z^2+z+1)t]\sigma_3},~~|x|\rightarrow \infty.
\end{equation}
Therefore, a transformation is given by
\begin{equation}
\Phi=\psi e^{-i[-zx+(z^3-z^2+z+1)t]\sigma_3},
\end{equation}
then the matrix function $\Phi=\Phi(x,t,z)$ has the following asymptotic property
\begin{equation}
\Phi \sim I,~~|x|\rightarrow \infty.
\end{equation}
and satisfies Lax pair
\begin{align}
&\Phi_x=-iz[\sigma_3,\Phi]+P\Phi,\label{lax2a}\\
&\Phi_t=i(z^3-z^2+z+1)[\sigma_3,\Phi]+Q\Phi.\label{lax2b}
\end{align}
The above two equations can be written in full differential form
\begin{equation}\label{eq1}
d\left(e^{i[zx-(z^3-z^2+z+1)t]\hat{\sigma}_{3}}\Phi \right)= e^{i[zx-(z^3-z^2+z+1)t]\hat{\sigma}_{3}}[Pdx+Qdt]\Phi.
\end{equation}
For asymptotic expansion of the eigenfunction of Lax pair \eqref{lax2a}-\eqref{lax2b} at infinity, we can proof that
\begin{align}
&\Phi \sim I, ~~z\rightarrow \infty,\\
&q(x,t)=2i \lim_{z\rightarrow \infty}(z\Phi)_{12}.
\end{align}
By integrating the equation \eqref{eq1} in two directions along the parallel real axis $(-\infty,t)\rightarrow (x,t)$ and $(+\infty,t)\rightarrow (x,t)$, we get two Volterra type integral equations
\begin{align}
&\Phi_{-}(x,t,z)=I+\int_{-\infty}^{x} e^{-iz(x-y)\hat{\sigma}_{3}}P(y,t,z)\Phi_{-}(y,t,z)dy ,\\
&\Phi_{+}(x,t,z)=I-\int_{x}^{+\infty} e^{-iz(x-y)\hat{\sigma}_{3}}P(y,t,z)\Phi_{+}(y,t,z)dy .
\end{align}
Due to $\Phi_{-}e^{-i[zx-(z^3-z^2+z+1)t]\sigma_{3}}$ and $\Phi_{+}e^{-i[zx-(z^3-z^2+z+1)t]\sigma_{3}}$
  being the linear correlation solutions of \eqref{lax1a} and \eqref{lax1b}, there exists a matrix $S(z)$ independent of $x$, $t$, which satisfies
\begin{equation}\label{s1}
\Phi_{-}(x,t,z)=\Phi_{+}e^{-i[zx-(z^3-z^2+z+1)t]\hat{\sigma}_{3}}S(z),
\end{equation}
where\\
$~~~~~~~~~~~~~~~~~~~~~~~~~~~~~~~~~~~~~~~~~~~~S(z)=\left(\begin{array}{cc}
s_{11}(z)  & s_{12}(z) \\
s_{21}(z)  &  s_{22}(z)
\end{array}\right)$.\\
Denote\\
$~~~~~~~~~~~~~~~~~~~~~~~~~~~~~~\Phi_-=(\Phi_{-1},\Phi_{-,2})$ and $\Phi_+=(\Phi_{+1},\Phi_{+,2})$\\
where the subscripts $1$ and $2$ represent the column vectors of the matrix respectively.\\
$\mathbf{Proposition~ 2.1.}$ $\Phi_{-,1}(z)$, $\Phi_{+,2}(z)$ and $s_{11}(z)$ are analytical in $\mathbb{C_+}$; $\Phi_{-,2}(z)$, $\Phi_{+,1}(z)$ and $s_{22}(z)$ are analytical in $\mathbb{C_-}$. Here, $\mathbb{C_+}=\{z\in \mathbb{C}|Imz>0\}$, $\mathbb{C_-}=\{z\in \mathbb{C}|Im z<0\}$.\\
$\mathbf{Proposition~ 2.2.}$ $\Phi_\pm(x,t,z)$ and $S(z)$ have the following symmetries
\begin{equation}\label{eq.2}
\Phi_{\pm}(x,t,z)=\sigma_{1} \overline{\Phi_{\pm}(x,t,\bar{z})} \sigma_{1}
\end{equation}
and
\begin{equation}\label{eq.3}
S(z)=\sigma_{1}\overline{S(\bar{z})}\sigma_{1},
\end{equation}
where\\
$~~~~~~~~~~~~~~~~~~~~~~~~~~~~~~~~~~~~~~~~~~~~~~~~~~~~\sigma_{1}=\left(\begin{array}{cc}
0  & 1 \\
1  & 0
\end{array}\right)$.\\
In addition, the scattering coefficients can be expressed in terms of the Jost functions as
\begin{equation}
s_{11}(z)=\det(\Phi_{-1},\Phi_{+2}),s_{22}(z)=\det(\Phi_{+1},\Phi_{-2}).
\end{equation}
Define reflection coefficient $r(z)=\frac{s_{21}(z)}{s_{11}(z)}$ and a sectionally meromorphic matrix
\begin{equation}
M(x,t,z)=\left\{\begin{array}{ll}
\left(\frac{\Phi_{-,1}(z)}{s_{11}(z)}, \Phi_{+,2}(z)\right), & \text {as} ~z \in \mathbb{C}^{+}, \\
\left(\Phi_{+,1}(z), \frac{\Phi_{-,2}(z)}{s_{22}(z)}\right), & \text {as} ~z \in \mathbb{C}^{-}.
\end{array}\right.
\end{equation}
The RH problem corresponding to the initial value problem of the NHNSE can be obtained from the equation \eqref{s1} by using the analytic property and symmetries of the characteristic function and spectral matrice.\\
\textbf{Riemann-Hilbert Problem 2.1:}\\
$\bullet$ Analyticity: $M(x,t,z)$ is analytical in $\mathbb{C} \setminus \mathbb{R}$;\\
$\bullet$ Jump condition: $M_{+}(x,t,z)=M_{-}(x,t,z)V(z)$, ~$z \in \mathbb{R}$; \\
$\bullet$ Asymptotic behaviors: $M(z) \rightarrow I$, ~$ z \rightarrow \infty$; \\
where\\
\begin{equation}\label{v}
V(z)=\left(\begin{array}{cc}
1-|r(z)|^{2} & -\overline{r}(z)e^{-2it\theta(z)} \\
r(z)e^{2it\theta(z)} & 1
\end{array}\right)
\end{equation}
and\\
$~~~~~~~~~~~~~~~~~~~~~~~~~~~~~~~~~~~~~~\theta(z)=z\frac{x}{t}-z^3+z^2-z-1 $.\\
This is an RH problem defined on the real axis, as shown in Figure 2.1. and the solution $q(x, t)$ of the initial value problem of the NHNSE can be described by the above RH problem.
\begin{equation}\label{q}
q(x,t)=2i\lim_{z\rightarrow \infty}(zM(z))_{12}=2i(M_1(z))_{12},
\end{equation}
where $M_1(z)$ comes from the asymptotic expansion of $M(z)$\\
$~~~~~~~~~~~~~~~~~~~~~~~~~~~~~~~~~~~~~~~~~M(z)=I+\frac{M_1}{z}+o(z^{-2}),~~z\rightarrow \infty $.\\

\centerline{\begin{tikzpicture}[scale=0.8]
\path  (1,0)--(-1,0) ;
\draw[-][thick](-4,0)--(-3,0);
\draw[-][thick](-4,0)--(-3,0);
\draw[-][thick](-3,0)--(-2,0);
\draw[-][thick](-2,0)--(-1,0);
\draw[-][thick](-1,0)--(0,0);
\draw[-][thick](-1,0)--(0,0);
\draw[-][thick](0,0)--(1,0);
\draw[-][thick](1,0)--(2,0);
\draw[-][thick](2,0)--(3,0);
\draw[-][thick](2,0)--(3,0);
\draw[->][thick](3,0)--(4,0)[thick]node[right]{$\mathbb{R}$};
\end{tikzpicture}}
~~~~~~~~~~~~~~~~~~~~~~~~~~~~~~~~~~~~\textbf{Figure 2.1.} Jump path of $M(z)$.

\section{Decompositions of the jump matrix $V(z)$}
We note that the long time asymptotic behavior of RH problem 2.1 is influenced by the decay and growth of the exponential function
\begin{equation}
e^{\pm2it\theta(z)},~~\theta(z)=z\frac{x}{t}-z^3+z^2-z-1.
\end{equation}
Based on this, we will analyse the real part of $\pm2it\theta(z)$ to ensure its exponential decay properties.
Therefore, in this section we introduce a new transformation $M(z) \rightarrow M^{(1)}(z)$ such that $M^{(1)}(z)$ behaves well along the characteristic line (or steepest line) as $t \rightarrow \infty$. Let $\xi=\frac{x}{t}$, to obtain the asymptotic properties of $e^{2it\theta(z)}$, we give the real part of $2it\theta(z)$:
\begin{align}\label{z}
&Im\theta(z)=\xi Imz-Im(z^3)+Im(z^2)-Imz,\\
&Re(2it\theta(z))=-2t\xi Imz+2tIm(z^3)-2tIm(z^2)+2tImz.
\end{align}

\begin{figure}[H]
\vspace{0cm}  
\centering
\subfloat[]{%
\begin{minipage}[]{0.3\textwidth}
\includegraphics[height=2.0in, width=2.0in]{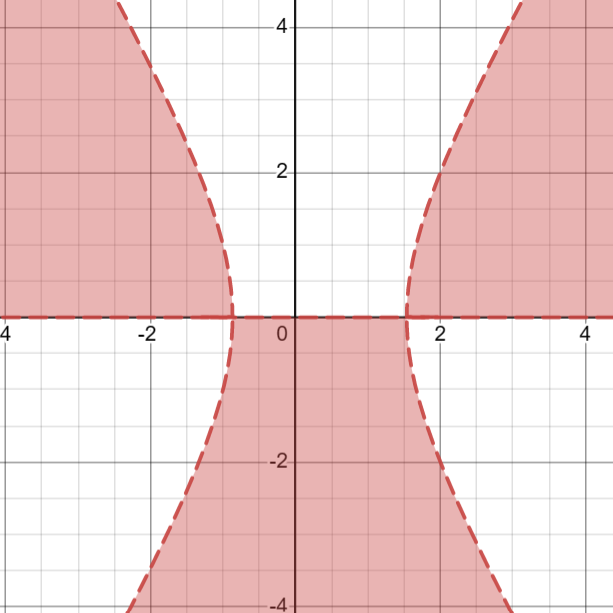}
\end{minipage}}
\caption*{Figure.3.1 The classification of sign Re$(i\theta(z))$.  In the white regions, Re$(i\theta(z)) < 0$, so $|e^{2it\theta(z)}| \rightarrow 0$ as $t \rightarrow \infty$. In the pink regions, Re$(i\theta(z)) > 0$, $|e^{-2it\theta(z)}| \rightarrow 0$ as $t \rightarrow \infty$.}
\end{figure}

To find the stationary phase points, we need the $\frac{d\theta}{dz}$
\begin{equation}
\frac{d\theta}{dz}=\xi-3z^2+2z-1.
\end{equation}
We get the stationary phase points
\begin{equation}
z_1=\frac{1+\sqrt{1+3(\xi-1)}}{3},~~z_2=\frac{1-\sqrt{1+3(\xi-1)}}{3}.
\end{equation}
When $\xi>\frac{2}{3}$, there are two stationary phase points; when $\xi<\frac{2}{3}$, there is no stationary phase point. In this paper, we consider the long-time asymptotics of a higher-order nonlinear Schr\"{o}dinger equation when $\xi>\frac{2}{3}$.\\
The jump matrix \eqref{v} has two decompositions
\begin{equation}\label{v01}
v(x,t,z)=\left(\begin{array}{cc}
1 & -\bar{r} e^{-2it \theta} \\
0 & 1
\end{array}\right)\left(\begin{array}{cc}
1 & 0 \\
r e^{2it \theta} & 1
\end{array}\right), z\in (z_{2},z_1),
\end{equation}
\begin{equation}\label{v02}
v(x,t,z)=\left(\begin{array}{cc}
1 & 0 \\
\frac{r}{1-|r|^2} e^{2it \theta} & 1
\end{array}\right)
\left(\begin{array}{cc}
1-|r|^2 & 0 \\
0  & \frac{1}{1-|r|^2}
\end{array}\right)
\left(\begin{array}{cc}
1 & \frac{-\bar{r}}{1-|r|^2} e^{-2it \theta} \\
0 & 1
\end{array}\right), z\in (-\infty,z_{2})\cup(z_1,+\infty).
\end{equation}
To remove the intermediate matrix of the second decomposition, we introduce a scalar RH problem:\\
$~~~~~~~~~~~~~~~~~~~~~~~~~~~~\bullet$ $\delta(z)$ is analytical in $\mathbb{C} \setminus \mathbb{R}$,\\
$~~~~~~~~~~~~~~~~~~~~~~~~~~~~\bullet$ $\delta_{+}(z)=\delta_{-}(z)(1-|r|^2)$, ~ $z\in (-\infty,z_{2})\cup(z_1,+\infty)$, \\
$~~~~~~~~~~~~~~~~~~~~~~~~~~~~~~\delta_{+}(z)=\delta_{-}(z)$, ~$z\in (z_{2},z_1)$,\\
$~~~~~~~~~~~~~~~~~~~~~~~~~~~~\bullet$ $\delta(z)\rightarrow 1$,~$z\rightarrow \infty$.\\
Using the Plemelj formula, it can be shown that this RH problem has a unique solution
\begin{equation}
\delta(z)=\exp \left[\frac{1}{2 \pi i} \int_{\Upsilon} \frac{\ln \left(1-|r(s)|^{2}\right)}{s-z} ds\right]=\exp \left[i \int_{\Upsilon} \frac{\nu(s)}{s-z} ds\right],
\end{equation}
where \\
$~~~~~~~~~~~~~~~~~~\Upsilon=(-\infty,z_{2})\cup(z_1,+\infty)$ and $\nu(s)=-\frac{1}{2\pi}\ln(1-|r(s)|^{2})$.\\
Suppose that $r(z)\in L^\infty \cap L^2$, and $\|r(z)\|_{L^\infty}\leq \rho <1$, then $\delta(z)$ has the following properties:\\
(1)$\delta(z)\overline{\delta(\bar{z})}=1$, $\|\delta_{\pm}-1\|_{L^2}\leq \frac{c\|r\|_{L^2}}{1-\rho}$;\\
(2)$(1-\rho^2)^{\frac{1}{2}}\leq |\delta(z)|\leq (1-\rho^2)^{-\frac{1}{2}}$.\\
 $\delta(z)$ can be rewritten as
\begin{equation}
\begin{aligned}\label{delta}
\delta(z)&=\exp\left(i\int_{\Upsilon}\frac{\nu(s)}{s-z}ds \right)\\
&=\exp\left(i\int_{\Upsilon}\frac{\nu(s)-\chi(s)\nu(z_j)}{s-z}ds+i\int_{z_{j}-1}^{z_{j}}\frac{\nu(z_j)}{s-z}ds\right)\\
&=(z-z_{j})^{i\nu(z_j)}\exp\left(i\int_{\Upsilon}\frac{\nu(s)-\chi(s)\nu(z_j)}{s-z}ds-i\nu(z_j)\log(z-z_j+1)\right)\\
&=(z-z_{j})^{i\nu(z_j)}\exp\left(i\lambda(z,z_j)\right),
\end{aligned}
\end{equation}
where $j=1,2$, $\chi(s)$ is the characteristic function defined on $(z_{j}-1, z_{j})$, \\
and\\
 $~~~~~~~~~~~~~~~~~~~z-z_{j}+1=|z-z_{j}+1|e^{i\arg(z-z_{j}+1)}$,\\
$~~~~~~~~~\exp(i\lambda(z,z_j))=\exp\left(i\int_{\Upsilon}\frac{\nu(s)-\chi(s)\nu(z_j)}{s-z}ds-i\nu(z_j)\log(z-z_j+1)\right)$.\\
We can prove that
\begin{align}
&\|\lambda(z,z_j)\|_{L^\infty}\leq \frac{c\|r\|_{H^{1,0}}}{1-\rho},\\
&|\lambda(z,z_j)-\lambda(z_j,z_{j})|<\frac{c\|r\|_{H^{1,0}}}{1-\rho}|z-z_{j}|^{\frac{1}{2}}.
\end{align}
Make a transformation
\begin{equation}\label{t1}
M^{(1)}(z)=M(z)\delta(z)^{-\sigma_3},
\end{equation}
$M^{(1)}(z)$ satisfies the following RH problem:\\
\textbf{Riemann-Hilbert Problem 3.1:}\\
$\bullet$ Analyticity: $M^{(1)}(x,t,z)$ is analytical in $\mathbb{C} \setminus \Sigma^{(1)}$;\\
$\bullet$ Jump condition: $M^{(1)}_{+}(x,t,z)=M^{(1)}_{-}(x,t,z)V^{(1)}(z), \quad z \in \Sigma^{(1)}$; \\
$\bullet$ Asymptotic behaviors: $M^{(1)}(z) \rightarrow I$, ~$ z \rightarrow \infty$; \\
where\\
$~~~~~~~~~~~~~~~~~~~~~~~~~~~~~~~~~~~~~~~~~~~~~~~~~~~~~~~~~~\Sigma^{(1)}=\mathbb{R}$,\\
and\\
\begin{equation}\label{v11}
V^{(1)}(z)=\left(\begin{array}{cc}
1 & -\delta(z)^2\bar{r} e^{-2it\theta(z)} \\
0 & 1
\end{array}\right)
\left(\begin{array}{cc}
1 & 0 \\
\delta(z)^{-2}re^{2it\theta(z)} & 1
\end{array}\right),~z\in(z_2,z_1);
\end{equation}

\begin{equation}\label{v12}
V^{(1)}(z)=\left(\begin{array}{cc}
1 & 0 \\
\delta_-(z)^{-2}\frac{r}{1-|r|^2}e^{2it\theta(z)} & 1
\end{array}\right)
\left(\begin{array}{cc}
1 & -\delta_+(z)^{2}\frac{\bar{r}}{1-|r|^2}e^{-2it\theta(z)}  \\
0 & 1
\end{array}\right),~z\in(-\infty,z_2)\cup (z_1,+\infty).
\end{equation}
Therefore, the relation between the solution of the NHNSE and the solution of the RH problem presents
\begin{equation}
q(x,t)=2i\lim_{z\rightarrow \infty}(zM^{(1)}(z)\delta(z)^{\sigma_3})_{12}=2i\lim_{z\rightarrow \infty}(zM^{(1)}(z))_{12}.
\end{equation}
\centerline{\begin{tikzpicture}[scale=0.8]
\path  (1,0)--(-1,0) ;
\draw[-][thick](-6,0)--(-5,0);
\draw[->][thick](-5,0)--(-4,0);
\draw[-][thick](-4,0)--(-3,0);
\draw[-][thick](-3,0)--(-2,0)node[below left]{$z_{2}$};
\draw[fill] (-2.5,0) circle [radius=0.05];
\draw[-][thick](-2,0)--(-1,0);
\draw[->][thick](-1,0)--(0,0);
\draw[-][thick](0,0)--(1,0);
\draw[-][thick](1,0)--(2,0)node[below left]{$z_{1}$};
\draw[fill] (1.6,0) circle [radius=0.05];
\draw[-][thick](2,0)--(3,0);
\draw[->][thick](3,0)--(4,0);
\draw[-][thick](4,0)--(6,0)[thick]node[right]{$\Sigma^{(1)}$};
\end{tikzpicture}}
~~~~~~~~~~~~~~~~~~~~~~~~~~~~~~~~~~~~~~~~~~~\textbf{Figure 3.2.} The oriented contour ${\Sigma^{(1)}}$.

\section{Continuous extension of scattering data}
We will deform jump contour of RH problem 3.1 to make the oscillating term of jump matrix being bounded on the transformed contour. Define path\\
$~~~~~~~~~~~~~~~~~~~~~~~~\Sigma^{(2)}=\Sigma_{11}\cup\Sigma_{12}\cup\Sigma_{13}\cup\Sigma_{14}\cup\Sigma_{21}\cup\Sigma_{22}\cup\Sigma_{23}
\cup\Sigma_{24}$,\\
where\\
$~~~~~~\Sigma_{11}=z_1+e^{i\varphi}R_{+}$, $\Sigma_{12}=z_1+e^{i(\pi-\varphi)}R_{+}$, $\Sigma_{13}=z_1+e^{i(\pi+\varphi)}R_{+}$, $\Sigma_{14}=z_1+e^{-i\varphi}R_{+}$,\\
$~~~~~~\Sigma_{21}=z_2+e^{i(\pi-\varphi)}R_{+}$, $\Sigma_{22}=z_2+e^{i\varphi}R_{+}$,
$\Sigma_{23}=z_2+e^{-i\varphi}R_{+}$, $\Sigma_{24}=z_2+e^{i(\pi+\varphi)}R_{+}$.\\

\centerline{\begin{tikzpicture}[scale=0.8]
\draw[->][dashed][thick](-8,0)--(8,0);
\draw[fill] (8.5,0) node{$Rez$};
\draw[fill] [blue](-3,0) circle [radius=0.05];
\draw[fill][blue] (-3,-0.3) node{$z_{2}$};
\draw[->][thick](-6,1.5)--(-4,0.5);
\draw[-][thick](-4,0.5)--(-3,0);
\draw[fill] (-4.5,0.3) node{$\Omega_{21}$};
\draw[fill] (-3,0.5) node{$\Omega_{22}$};
\draw[fill] (-1.5,0.3) node{$\Omega_{23}$};
\draw[fill] (-1.5,-0.3) node{$\Omega_{24}$};
\draw[fill] (-3,-0.8) node{$\Omega_{25}$};
\draw[fill] (-4.5,-0.3) node{$\Omega_{26}$};
\draw[fill] (-6.5,1.7) node{$\Sigma_{21}$};
\draw[fill] (-6.5,-1.7) node{$\Sigma_{24}$};
\draw[fill] (-1.5,1.2) node{$\Sigma_{22}$};
\draw[fill] (-1.5,-1.2) node{$\Sigma_{23}$};
\draw[->][thick](-6,-1.5)--(-4,-0.5);
\draw[-][thick](-4,-0.5)--(-3,0);
\draw[->][thick](-3,0)--(-1.5,0.75);
\draw[-][thick](-1.5,0.75)--(-0,1.5);
\draw[->][thick](-3,0)--(-1.5,-0.75);
\draw[-][thick](-1.5,-0.75)--(0,-1.5);
\draw[fill][red] (0,0) circle [radius=0.05];
\draw[fill][red] (0,-0.3) node{$0$};
\draw[->][dashed][thick](0,-3)--(0,3);
\draw[fill] (0.6,3) node{$Imz$};
\draw[fill] [blue](-3,0) circle [radius=0.05];
\draw[fill][blue] (-3,-0.3) node{$z_{2}$};
\draw[-][thick](6,1.5)--(4,0.5);
\draw[<-][thick](4,0.5)--(3,0);
\draw[fill] (4.5,0.3) node{$\Omega_{11}$};
\draw[fill] (3,0.5) node{$\Omega_{12}$};
\draw[fill] (1.5,0.3) node{$\Omega_{13}$};
\draw[fill] (1.5,-0.3) node{$\Omega_{14}$};
\draw[fill] (3,-0.8) node{$\Omega_{15}$};
\draw[fill] (4.5,-0.3) node{$\Omega_{16}$};
\draw[fill] (6.5,1.7) node{$\Sigma_{11}$};
\draw[fill] (6.5,-1.7) node{$\Sigma_{14}$};
\draw[fill] (1.5,1.2) node{$\Sigma_{12}$};
\draw[fill] (1.5,-1.2) node{$\Sigma_{13}$};
\draw[-][thick](6,-1.5)--(4,-0.5);
\draw[<-][thick](4,-0.5)--(3,0);
\draw[-][thick](3,0)--(1.5,0.75);
\draw[<-][thick](1.5,0.75)--(0,1.5);
\draw[-][thick](3,0)--(1.5,-0.75);
\draw[<-][thick](1.5,-0.75)--(0,-1.5);
\draw[fill][blue] (3,-0.3) node{$z_{1}$};
\end{tikzpicture}}
\textbf{Figure 4.1.} Deformation from $\mathbb{R}$ to $\Sigma^{(2)}$.

We introduce another matrix valued function $M^{(2)}(z)$:
\begin{equation}\label{t2}
M^{(2)}(z)=M^{(1)}(z)R^{(2)}(z).
\end{equation}
Here\\
\begin{equation}\label{r2}
R^{(2)}(z)=\begin{cases}
\left(
\begin{array}{cc}
1 & -R_{j1}(z)\delta(z)^2 e^{-2it \theta} \\
0 & 1 \\
 \end{array}
 \right) , ~~as~z \in \Omega_{j1};\\
 \left(
\begin{array}{cc}
1 & 0 \\
-R_{j3}(z)\delta(z)^{-2} e^{2it\theta} & 1 \\
 \end{array}
 \right) , ~~as~z \in \Omega_{j3};\\
 \left(
\begin{array}{cc}
1 & R_{j4}(z)\delta(z)^2e^{-2it \theta} \\
0 & 1 \\
 \end{array}
 \right) , ~~as ~z \in \Omega_{j4};\\
 \left(
\begin{array}{cc}
1 & 0 \\
R_{j6}(z)\delta(z)^{-2} e^{-2it \theta} & 1 \\
 \end{array}
 \right) , ~~as~z \in \Omega_{j6};\\
 I,~~~~~~~~~~~~~~~~~~~~~~~~~~~~~~~~~~~~~as~z~in~elsewhere;\\
\end{cases}
\end{equation}
where $j=1,2$.\\
\textbf{Proposition~4.1} There exist functions $R_{jk}: \Omega_{jk} \rightarrow \mathbb{C}$, $j=1,2$, $k=1,3,4,6$ satisfying the following boundary conditions
\begin{equation}\label{Rj1}
R_{j1}(z)=\begin{cases}
-\frac{\bar{r}(z)}{1-|r(z)|^2}, ~z \in (-\infty,z_2)\cup (z_1,+\infty);\\
 f_{j1}=-\frac{\bar{r}(z_j)}{1-|r(z_j)|^2}e^{2i\lambda(z_j,z_j)}(z-z_j)^{2i\nu(z_j)}\delta(z)^{-2} , ~z \in \Sigma_{j1};\\
\end{cases}
\end{equation}

\begin{equation}\label{Rj3}
R_{j3}(z)=\begin{cases}
r(z), ~z \in (z_2,z_1);\\
 f_{j3}=r(z_j)e^{-2i\lambda(z_j,z_j)}(z-z_j)^{-2i\nu(z_j)}\delta(z)^{2} , ~z \in \Sigma_{j2};\\
\end{cases}
\end{equation}

\begin{equation}\label{Rj4}
R_{j4}(z)=\begin{cases}
-\bar{r}(z), ~z \in (z_2,z_1);\\
 f_{j4}=-\bar{r}(z_j)e^{2i\lambda(z_j,z_j)}(z-z_j)^{2i\nu(z_j)}\delta(z)^{-2} , ~z \in \Sigma_{j3};\\
\end{cases}
\end{equation}

\begin{equation}\label{Rj6}
R_{j6}(z)=\begin{cases}
\frac{r(z)}{1-|r(z)|^2}, ~z \in (-\infty,z_2)\cup (z_1,+\infty);\\
 f_{j6}=\frac{r(z_j)}{1-|r(z_j)|^2}e^{-2i\lambda(z_j,z_j)}(z-z_j)^{-2i\nu(z_j)}\delta(z)^{2} , ~z \in \Sigma_{j4}.\\
\end{cases}
\end{equation}
$R_{jk}(z)$ has the following estimations
\begin{align}
&|R_{jk}(z)|\leq c_1\sin^{2}(\arg(z-z_j))+c_2\langle Rez\rangle ^{-1},\\
&\bar{\partial}R_{jk}(z)\leq|f_{jk}^{'}(Rez)|+|z-z_j|^{-\frac{1}{2}},
\end{align}
where $\langle Rez\rangle=\left(1+(Rez)^2  \right)^{\frac{1}{2}}$.\\
$\mathbf{Proof.}$ Taking the $R_{11}(z)$ as an example to show this proof.\\
Define\\
$~~~~~~~~~~~~~~~~~~~~~~~~~~~~~~~~~~~~z=z_1+\rho e^{i\varphi}$, $\bar{z}=z_1+\rho e^{-i\varphi}$, \\
the derivative of the above two equations with respect to $\bar{z}$ gives
\begin{equation}
\frac{\partial \rho}{\partial \bar{z}}+i\rho \frac{\partial \varphi}{\partial \bar{z}}=0,
\frac{\partial \rho}{\partial \bar{z}}-i\rho \frac{\partial \varphi}{\partial \bar{z}}=e^{i\varphi},
\end{equation}
the resulting solution is
\begin{equation}
\frac{\partial \rho}{\partial \bar{z}}=\frac{1}{2}e^{i\varphi},
 \frac{\partial \varphi}{\partial \bar{z}}=\frac{ie^{i\varphi}}{2\rho}.
\end{equation}
Thus, the full derivative of $\bar{\partial}$ in polar coordinates is
\begin{equation}\label{polar}
\bar{\partial}=\frac{1}{2}e^{i\varphi}(\partial_{\rho}+i\rho^{-1}\partial_{\varphi}).
\end{equation}
Defining function
\begin{equation}\label{R11}
R_{11}(z)=\cos(2\varphi)p_{11}(Rez)+[1-\cos(2\varphi)]f_{11}(z),
\end{equation}
where $z=z_{1}+\rho e^{i\varphi}$, $\rho \geq 0$, $0 \leq \varphi \leq \frac{\pi}{4}$, $p_{11}(z)=\frac{\bar{r}(z)}{1-|r(z)|^2}$, then\\
$~~~~~~~~~~~~~~~~~~~~~~~~~~~~~~~~~\rho=|z-z_1|$, $Rez=z_1+\rho \cos\varphi$.\\
For $z\in (z_1,+\infty)$, $\varphi=0$, thus $R_{11}(z)=p_{11}(Rez)$; for $z\in \Sigma_{11}$, $\varphi=\frac{\pi}{4}$, so $R_{11}(z)=f_{11}(z)$.\\
Notice that $f_{11}(z)$ is analytical, using \eqref{polar}, then we have
\begin{equation}
\begin{aligned}
\bar{\partial}R_{11}(z)&=(p_{11}(Rez)-f_{11}(z))\bar{\partial}\cos(2\varphi)+\frac{1}{2}e^{i\varphi}\cos(2\varphi)
p_{11}^{'}(Rez)\\
&=i(p_{11}(Re(z))-f_{11}(z))\frac{1}{\rho}e^{i\varphi}\sin(2\varphi)+\frac{1}{2}e^{i\varphi}\cos(2\varphi)p_{11}^{'}(Rez),
\end{aligned}
\end{equation}
therefore,
\begin{equation}
|\bar{\partial}R_{11}(z)|\leq \frac{c_1}{|z-z_1|}(|p_{11}(z)-p_{11}(z_1)|+|p_{11}(z_1)-f_{11}(z)|)+c_2|p_{11}^{'}(Rez)|.
\end{equation}
While
\begin{equation}
\begin{aligned}
|p_{11}(z)-p_{11}(z_1)|&=|\int_{z_1}^{z}p_{11}^{'}(s)ds|\leq \int_{z_1}^{z}|p_{11}^{'}(s)|ds\\
&\leq \|p_{11}^{'}(s) \|_{L^2(z,z_1)} \|1\|_{L^2(z,z_1)}\\
&\leq c|z-z_1|^{\frac{1}{2}}.\\
\end{aligned}
\end{equation}
Using \eqref{delta}, we get
\begin{equation}
\begin{aligned}
f_{11}(z)&=p_{11}(z_1)e^{2i\lambda(z_1,z_1)}(z-z_1)^{2i\nu(z_1)}\delta(z)^{-2}\\
&=p_{11}(z_1)e^{2i\lambda(z_1,z_1)}e^{-2i\lambda(z,z_1)},
\end{aligned}
\end{equation}
then
\begin{equation}
p_{11}(z_1)-f_{11}(z)=p_{11}(z_1)-p_{11}(z_1)e^{2i\lambda(z_1,z_1)}e^{-2i\lambda(z,z_1)},
\end{equation}
 $z \in \Omega_{11}$, it holds that
\begin{equation}
|\lambda(z,z_1)-\lambda(z_1,z_1)|=O(|z-z_1|^{\frac{1}{2}}),
\end{equation}
hence
\begin{equation}
|p_{11}(z_1)-f_{11}(z)|=p_{11}(z_1)(1-\exp[O(|z-z_1|^{\frac{1}{2}})])=p_{11}(z_1)O(|z-z_1|^{\frac{1}{2}}).
\end{equation}
Thus, we obtain that
\begin{equation}
\bar{\partial}R_{11}(z)\leq c_1 |z-z_1|^{-\frac{1}{2}}+c_2|p_{11}^{'}(Rez)|.
\end{equation}
From \eqref{R11}, it is easy to have that
\begin{equation}
\begin{aligned}
|R_{11}|&\leq 2|f_{11}(z)|\sin^2(\varphi)+|\cos(2\varphi)||p_{11}(Rez)|\\
&\leq c_1\sin^2(\varphi)+c_2[1+(Rez)^2]^{-\frac{1}{2}}\\
&=c_1\sin^2(\varphi)+c_2\langle Rez \rangle^{-1}.
\end{aligned}
\end{equation}
$M^{(2)}(x,t,z)$ satisfies the following RH problem:\\
\textbf{Riemann-Hilbert Problem 4.1:}\\
$\bullet$ Analyticity: $M^{(2)}(x,t,z)$ is analytical in $\mathbb{C} \setminus \Sigma^{(2)}$;\\
$\bullet$ Jump condition: $M^{(2)}_{+}(x,t,z)=M^{(2)}_{-}(x,t,z)V^{(2)}(z), ~ z \in \Sigma^{(2)}$; \\
$\bullet$ Asymptotic behaviors: $M^{(2)}(z) \rightarrow I$, ~$ z \rightarrow \infty$ ;\\
where\\
\begin{equation}\label{V2}
V^{(2)}(z)=\begin{cases}
\left(
\begin{array}{cc}
1 & R_{j1}(z)\delta(z)^2 e^{-2it \theta} \\
0 & 1 \\
 \end{array}
 \right) , ~z \in \Sigma_{j1};\\
 \left(
\begin{array}{cc}
1 & 0 \\
R_{j3}(z)\delta(z)^{-2} e^{2it\theta} & 1 \\
 \end{array}
 \right) , ~z \in \Sigma_{j2};\\
 \left(
\begin{array}{cc}
1 & R_{j4}(z)\delta(z)^2e^{-2it \theta} \\
0 & 1 \\
 \end{array}
 \right) ,  ~z \in \Sigma_{j3};\\
 \left(
\begin{array}{cc}
1 & 0 \\
R_{j6}(z)\delta(z)^{-2} e^{2it \theta} & 1 \\
 \end{array}
 \right) ,~z \in \Sigma_{j4};\\
 I,~~~~~~~~~~~~~~~~~~~~~~~~~~~~~~~~~z\in (z_2,z_1)\cup (z_1,+\infty),\\
\end{cases}
\end{equation}
here $j=1,2$.\\
The relation between the solution of the higher-order nonlinear Schr\"{o}dinger equation and the solution of the RH problem is that
\begin{equation}
q(x,t)=2i\lim_{z\rightarrow \infty}(zM^{(2)}(z))_{12}.
\end{equation}

\section{$\bar{\partial}$-problems and solution asymptotics}
Denote\\
\begin{equation}\label{t3}
M^{(3)}(z)=M^{(2)}(z)M^{(lo)}(z)^{-1}.
\end{equation}
$M^{(lo)}(z)$ satisfies the following RH problem:\\
\textbf{Riemann-Hilbert Problem 5.1:}\\
$\bullet$ Analyticity: $M^{(lo)}(x,t,z)$ is analytical in $\mathbb{C} \setminus \Sigma^{(0)}$;\\
$\bullet$ Jump condition: $M^{(lo)}_{+}(x,t,z)=M^{(lo)}_{-}(x,t,z)V^{(2)}(z), ~ z \in \Sigma^{(0)}$; \\
$\bullet$ Asymptotic behaviors: $M^{(lo)}(z) \rightarrow I$, ~$ z \rightarrow \infty$;\\
where $\Sigma^{(0)}=(\mathop{\cup}\limits_{j=1,2;k=1,2,3,4}\Sigma_{jk})\cap U(z_j)$. The jump matrix $V^{(2)}(z)$ has the following decomposition
\begin{equation}
V^{(2)}(z)=(I-\omega_{jk}^{-}(z))^{-1}(I+\omega_{jk}^{+}(z)),
\end{equation}
where $\omega_{jk}^{+}(z)=0$, $\omega_{jk}^{-}(z)=I-V^{(2)}(z)^{-1}=V^{(2)}(z)-I$,
\begin{equation}
w_{jk}^{-}(z)=\left\{\begin{array}{cc}
\left(\begin{array}{cc}
0 & R_{jk}(z)\delta(z)^2 e^{-2it \theta} \\
0 & 0
\end{array}\right), & z \in \Sigma_{jk},j=1,2;k=1,3,  \\
\left(\begin{array}{cc}
0 & 0 \\
R_{jk}(z)\delta(z)^{-2} e^{2it \theta} & 0
\end{array}\right), & z \in \Sigma_{jk},j=1,2;k=2,4, \\
\end{array}\right.
\end{equation}
Define\\
$~~~~~~~~~~~~~~~~~~~~~~~~~\Sigma_j^{(0)}=\mathop{\cup}\limits_{k=1,2,3,4}\Sigma_{jk})$, $\omega_{j}(z)=\mathop{\Sigma}\limits_{k=1,2,3,4}\omega_{jk}(z)$, \\
$~~~~~~~~~~~~~~~~\omega_{jk}^{\pm}(z)=\omega_{jk}(z)|_{C^{\pm}}$,
$\omega_{k}^{\pm}(z)=\omega_{k}(z)|_{C^{\pm}}$,
$\omega^{\pm}(z)=\omega(z)|_{C^{\pm}}$.\\
The Cauchy projection operator $C_{\pm}$ on $\Sigma^{(2)}$ is
\begin{equation}
C_{\pm}(f)(s)=\lim _{z \rightarrow \Sigma_{\pm}^{(2)}} \frac{1}{2 \pi i} \int_{\Sigma^{(2)}} \frac{f(s)}{s-z} d s,
\end{equation}
based on the above, we define operator
\begin{align}
&C_{\omega}(f)=C_{+}\left(f\omega^{-}\right)+C_{-}\left(f\omega^{+}\right),
C_{\omega_j}(f)=C_{+}\left(f\omega_{j}^{-}\right)+C_{-}\left(f\omega_{j}^{+}\right),\\
&C_{\omega}=\Sigma_{j=1}^2C_{\omega_j}.
\end{align}
$\mathbf{Lemma~ 5.1.}$ The matrix functions $\omega_{jk}^{-}$ defined above admits following estimation:
\begin{equation}
\left\|\omega_{jk}^{-}\right\|_{L^{p}\left(\Sigma_{jk}\right)}=\mathcal{O}\left(t^{-1/2}\right), 1 \leq p<+\infty.
\end{equation}
From Lemma 5.1, we know that $I-C_\omega$ and $I-C_{\omega_j}$ are reversible, then the solution of RH problem 5.1 is unique. By using Beals-Cofiman theorem, the solution can be written as
\begin{equation}
M^{(lo)}(z)=I+\frac{1}{2 \pi i} \int_{\Sigma^{(0)}} \frac{\left(I-C_{w}\right)^{-1} I \cdot w}{s-z} ds.
\end{equation}
$\mathbf{Corollary~5.1}$ As $t\rightarrow \infty$,
\begin{equation}
\|C_{\omega_{j}}C_{\omega_{k}}\|_{L^2(\Sigma^{(0)})}\leq t^{-1},
\|C_{\omega_{j}}C_{\omega_{k}}\|_{L^{\infty}(\Sigma^{(0)})\rightarrow L^2(\Sigma^{(0)})}\leq t^{-1}.
\end{equation}
By direct calculation, we have\\
$\left(I-C_{w}\right)\left(I+\sum_{j=1}^{2} C_{w_{j}}\left(I-C_{w_{j}}\right)^{-1}\right)=I-\sum_{1 \leq j \neq k \leq 2} C_{w_{k}} C_{w_{j}}\left(I-C_{w_{j}}\right)^{-1} $,\\
$\left(I+\sum_{j=1}^{2} C_{w_{j}}\left(I-C_{w_{j}}\right)^{-1}\right)\left(I-C_{w}\right)=I-\sum_{1 \leq j \neq k \leq 2}\left(I-C_{w_{j}}\right)^{-1} C_{w_{j}} C_{w_{k}}$.\\
$\mathbf{Proposition ~5.1}$ As $t\rightarrow \infty$, we obtain
\begin{equation}
\int_{\Sigma^{(0)}} \frac{\left(I-C_{\omega}\right)^{-1} I \cdot \omega}{s-z} d s=\sum_{j=1}^{2} \int_{\Sigma_{j}^{(0)}} \frac{\left(I-C_{\omega_{j}}\right)^{-1} I \cdot \omega_{j}}{s-z} ds+\mathcal{O}\left(t^{-3 / 2}\right).
\end{equation}
Thus, as $t\rightarrow \infty$, we consider the reduction of RH problem 5.1 to a model RH problem whose solutions can be given explicitly in terms of parabolic cylinder functions on each contour $\Sigma_j^{(0)}$ respectively. We denote $\hat{\Sigma}_{1}^{(0)}$ as the contour $z=z_1+le^{\pm i\varphi},l\in R$ oriented from $\Sigma_1^{(0)}$, and $\hat{\Sigma}_{1k}$ is the extension of $\Sigma_{1k}$ respectively. And for $z$ near $z_1$, rewrite phase
function as
\begin{equation}
\theta(z)=\theta\left(z_{1}\right)+\left(z-z_{1}\right)^{2} \theta^{\prime \prime}\left(z_{1}\right)+\mathcal{O}\left(\left(z-z_{1}\right)^{3}\right).
\end{equation}
\textbf{Riemann-Hilbert Problem 5.2:} Find a matrix-valued function $M^{(lo,1)}(x,t,z)$ with following properties\\
$\bullet$ Analyticity: $M^{(lo,1)}(x,t,z)$ is analytical in $\mathbb{C} \setminus \hat{\Sigma}_{1}^{(0)}$;\\
$\bullet$ Jump condition: $M^{(lo,1)}_{+}(x,t,z)=M^{(lo,1)}_{-}(x,t,z)V^{lo,1}(z), ~ z \in \hat{\Sigma}_{1}^{(0)}$; \\
$\bullet$ Asymptotic behaviors: $M^{(lo,1)}(z) \rightarrow I$, ~$ z \rightarrow \infty$;\\
where
\begin{equation}
V^{lo,1}(z)=\left\{\begin{array}{cc}
\left(\begin{array}{cc}
1 & -\frac{\bar{r}(z_1)}{1-|r(z_1)|^2}e^{2i\lambda(z_1,z_1)} (z-z_1)^{2i\nu(z_1)}e^{-2it \theta} \\
0 & 1
\end{array}\right), & z \in \hat{\Sigma}_{11},  \\
\left(\begin{array}{cc}
1 & 0 \\
r(z_1)e^{-2i\lambda(z_1,z_1)} (z-z_1)^{-2i\nu(z_1)}e^{2it \theta} & 1
\end{array}\right), & z \in \hat{\Sigma}_{12},  \\
\left(\begin{array}{cc}
1 & -\bar{r}(z_1)e^{2i\lambda(z_1,z_1)} (z-z_1)^{2i\nu(z_1)}e^{-2it \theta} \\
0 & 1
\end{array}\right), & z \in \hat{\Sigma}_{13},  \\
\left(\begin{array}{cc}
1 & \frac{r(z_1)}{1-|r(z_1)|^2}e^{-2i\lambda(z_1,z_1)} (z-z_1)^{-2i\nu(z_1)}e^{2it \theta} \\
0 & 1
\end{array}\right), & z \in \hat{\Sigma}_{14},  \\
\end{array}\right.
\end{equation}

\centerline{\begin{tikzpicture}[scale=0.8]
\draw[dashed][thin](-6,0)--(6,0);
\draw[fill] (0,0) circle [radius=0.05];
\draw[fill] (0,-0.3) node{$z_{1}$};
\draw[->][thick](-4,4)--(-2,2);
\draw[-][thick](-2,2)--(0,0);
\draw[->][thick](-4,-4)--(-2,-2);
\draw[-][thick](-2,-2)--(0,0);
\draw[fill] (-4.5,2.3) node{$\hat{\Sigma}_{12}$};
\draw[fill] (-4.5,-2.3) node{$\hat{\Sigma}_{13}$};
\draw[fill] (-7.5,1) node{$\left(\begin{array}{cc}
1 & 0 \\
r(z_1)e^{-2i\lambda(z_1,z_1)}(z-z_1)^{-2i\nu (z_1)}e^{2it\theta} & 1
\end{array}\right)$};
\draw[fill] (-7.5,-1) node{$\left(\begin{array}{cc}
1 & -\bar{r}(z_1)e^{2i\lambda(z_1,z_1)}(z-z_1)^{2i\nu (z_1)}e^{-2it\theta} \\
0 & 1
\end{array}\right)$};
\draw[->][thick](0,0)--(2,2);
\draw[-][thick](2,2)--(4,4);
\draw[->][thick](0,0)--(2,-2);
\draw[-][thick](2,-2)--(4,-4);
\draw[fill] (4.5,2.3) node{$\hat{\Sigma}_{11}$};
\draw[fill] (4.5,-2.3) node{$\hat{\Sigma}_{14}$};
\draw[fill] (7.5,1) node{$\left(\begin{array}{cc}
1 & -\frac{\bar{r}(z_1)}{1+|r(z_1)|^2}e^{2i\lambda(z_1,z_1)}(z-z-1)^{2i\nu (z_1)}e^{-2it\theta} \\
0 & 1
\end{array}\right)$};
\draw[fill] (7.5,-1) node{$\left(\begin{array}{cc}
1 & 0 \\
\frac{r(z_1)}{1+|r(z_1)|^2}e^{-2i\lambda(z_1,z_1)}(z-z_1)^{-2i\nu (z_1)}e^{2it\theta} & 1
\end{array}\right)$};
\end{tikzpicture}}
\textbf{Figure 5.1.} The contour $\hat{\Sigma}_1^{(0)}$ and the jump matrix.\\
To match the model, the translational scale transformation is given by
\begin{equation}
\xi(z)=t^{1/2} \sqrt{2\theta^{\prime \prime}(z_{1})}\left(z-z_{1}\right).
\end{equation}
Let\\
$~~~~~~~~~~~~~~~~~~~~~~r_{z_{1}}=r\left(z_{1}\right)e^{-2i\lambda(z_1,z_1)+2it\theta(z_1)}\exp(-i\nu(z_1)\log(2t\theta^{''}(z_1)))$,\\
with $|r_{z_{1}}|=|r(z_{1})|$. Through this change of variable, the jump $V^{lo,1}(z)$ approximates to the jump of a
parabolic cylinder model problem as follows:\\
\textbf{Riemann-Hilbert Problem 5.3:} Find a matrix-valued function $M^{pc}(\xi)$ with following properties\\
$\bullet$ Analyticity: $M^{pc}(\xi)$ is analytical in $\mathbb{C} \setminus \Sigma^{pc}$;\\
$\bullet$ Jump condition: $M^{pc}_{+}(\xi)=M^{pc}_{-}(\xi)V^{pc}(\xi), ~ \xi \in \Sigma^{pc}$; \\
$\bullet$ Asymptotic behaviors: $M^{pc}(\xi)= I+\frac{M_{1}^{pc}}{\xi}+O(\xi^{-2})$, ~$\xi \rightarrow \infty$;\\
where
\begin{equation}
V^{pc}(\xi)=\left\{\begin{array}{cc}
\left(\begin{array}{cc}
1 & -\frac{\bar{r}_{z_1}}{1-|r_{z_1}|^2}\xi^{2i\nu}e^{-\frac{i}{2}\xi^2} \\
0 & 1
\end{array}\right), & \xi \in R^{+}e^{i\varphi},  \\
\left(\begin{array}{cc}
1 & 0 \\
r_{z_1}\xi^{-2i\nu}e^{\frac{i}{2}\xi^2} & 1
\end{array}\right), & \xi \in R^{+}e^{i(\pi-\varphi)},  \\
\left(\begin{array}{cc}
1 & -\bar{r}_{z_1}\xi^{2i\nu}e^{-\frac{i}{2}\xi^2} \\
0 & 1
\end{array}\right), & \xi \in R^{+}e^{i(\pi+\varphi)},  \\
\left(\begin{array}{cc}
1 & 0 \\
\frac{r_{z_1}}{1-|r_{z_1}|^2} \xi^{-2i\nu}e^{\frac{i}{2}\xi^2} & 1
\end{array}\right), & \xi \in R^{+}e^{-i\varphi},  \\
\end{array}\right.
\end{equation}
Make a transformation
\begin{equation}
M^{pc}(\xi)=\Psi(\xi)P(\xi)\xi^{-i\nu \sigma_3} e^{\frac{i}{4}\xi^{2} \sigma_{3}},
\end{equation}
where
\begin{equation}
P(\xi)=\left\{\begin{array}{ll}
\left(\begin{array}{cc}
1 & \frac{\bar{r}_{z_1}}{1-|r_{z_1}|^2} \\
0 & 1
\end{array}\right), & \xi \in \Omega_1, \\
\left(\begin{array}{cc}
1 & 0 \\
-r_{z_1} & 1
\end{array}\right), & \xi \in \Omega_3, \\
\left(\begin{array}{cc}
1 & -\bar{r}_{z_1} \\
0 & 1
\end{array}\right), & \xi \in \Omega_4, \\
\left(\begin{array}{cc}
1 & 0 \\
\frac{r_{z_1}}{1-|r_{z_1}|^2} & 1
\end{array}\right), & \xi \in \Omega_6, \\
I, & \xi \in \Omega_2 \cup \Omega_5,
\end{array}\right.
\end{equation}
then a standard RH problem with jumps only at $\xi = 0$ is obtained\\
\textbf{Riemann-Hilbert Problem 5.4:} Find a matrix-valued function $\Psi(\xi)$ which satisfies following properties:\\
$\bullet$ Analyticity: $\Psi(\xi)$ is analytical in $\mathbb{C}\setminus \mathbb{R}$;\\
$\bullet$ Jump condition: $\Psi(\xi)$ has continuous boundary values $\Psi_{\pm}(\xi)$ on $\mathbb{R}$ and
\begin{equation}\label{v0}
\Psi_{+}(\xi)=\Psi_{-}(\xi)V(0),~~\xi \in \mathbb{R},
\end{equation}
where\\
\begin{equation}
V(0)=\left(\begin{array}{cc}
1-|r_{z_1}|^2 & -\bar{r}_{z_{1}} \\
r_{z_{1}} & 1
\end{array}\right);
\end{equation}
$\bullet$ Asymptotic behaviors:
\begin{equation}
\Psi(\xi)e^{\frac{i}{4}\xi^2\sigma_3}\xi^{-i\nu\sigma_3}\rightarrow I, ~~\xi \rightarrow \infty.
\end{equation}
The above RH problem can be reduced to the Weber equation to obtain an explicit solution in the parabolic plane $\Psi(\xi)=(\Psi_{ij})_{i,j=1}^2$ which can be given by a parabolic function on the upper half-plane $\xi \in \mathbb{C}^+$.\\
Taking the derivative of both sides of \eqref{v0} with respect to $\xi$, we get
\begin{equation}
\left(\frac{d\Psi}{d\xi} \right)_{+}=\left(\frac{d\Psi}{d\xi} \right)_{-}V(0),
\end{equation}
then use\\
$~~~~~~~~~~~~~~~~~~~~~~~~~~~~~~~~~~~~~~~~~~~~~~~\frac{i}{2}\xi \sigma_3\Psi_+=\frac{i}{2}\xi \sigma_3\Psi_-V(0)$,\\
have
\begin{equation}
\left(\frac{d\Psi}{d\xi}+ \frac{i}{2}\xi \sigma_3\Psi\right)_{+}=\left(\frac{d\Psi}{d\xi}+ \frac{i}{2}\xi \sigma_3\Psi\right)_{-}V(0).
\end{equation}
Since $\det (V(0))=1$, taking the determinant on both sides of \eqref{v0} gives
\begin{equation}
\det(\Psi)_{+}=\det(\Psi)_{-}, ~\xi\in \mathbb{R}.
\end{equation}
By Painlev$\acute{e}$'s theorem, we know that $\det(\Psi)$ is analytical on $\mathbb{C}$ and $\det(\Psi)=\det(M^{pc}) \rightarrow 1$, so $\det(\Psi)$ is bounded, and $\det(\Psi)=1$ , then $\Psi$ is invertible, and
\begin{equation}\label{mpc-1}
(M^{pc})^{-1}=e^{-\frac{i}{4}\xi^2\sigma_3}\xi^{i\nu\sigma_3}P^{-1}\psi^{-1}=I-\frac{M_1^{pc}}{\xi}+O(\xi^{-2})
\end{equation}
is bounded. By using \eqref{mpc-1}, we have\\
$~~~~~~~~~~~~~~~~~~~~~~~~~~~~\left[(\frac{d\Psi}{d\xi}+ \frac{i}{2}\xi \sigma_3\Psi) \Psi^{-1}\right]_{+}=(\frac{d\Psi}{d\xi}+ \frac{i}{2}\xi \sigma_3\Psi)_{+}\Psi_{+}^{-1}$\\
$~~~~~~~~~~~~~~~~~~~~~~~~=(\frac{d\Psi}{d\xi}+ \frac{i}{2}\xi \sigma_3\Psi)_{-}V(0)V(0)^{-1}\Psi_{-}^{-1}=\left[(\frac{d\Psi}{d\xi}+ \frac{i}{2}\xi \sigma_3\Psi) \Psi^{-1}\right]_{-}$.\\
From Painlev$\acute{e}$'s theorem, we know that $(\frac{d\Psi}{d\xi}+ \frac{i}{2}\xi \sigma_3\Psi) \Psi^{-1}$ is analytical on $\mathbb{C}$. A direct calculation gives
\begin{equation}
\begin{aligned}
&~~\left(\frac{d\Psi}{d\xi}+ \frac{i}{2}\xi \sigma_3\Psi\right) \Psi^{-1}\\
&=\frac{dM^{pc}}{d\xi}(M^{pc})^{-1}-\frac{i}{2}\xi M^{pc}\sigma_3 (M^{pc})^{-1}+i\xi M^{pc}\sigma_3 (M^{pc})^{-1}
+\frac{i}{2}\xi \sigma_3 M^{pc} (M^{pc})^{-1}\\
&=\frac{i}{2}\xi[\sigma_3, M^{pc}](M^{pc})^{-1}+O(\xi^{-1})\\
&=\frac{i}{2}[\sigma_3, M_{1}^{pc}]+O(\xi^{-1})
\end{aligned}
\end{equation}
is bounded. By Liouville's theorem, $\left(\frac{d\Psi}{d\xi}+ \frac{i}{2}\xi \sigma_3\Psi\right)\Psi^{-1}$ is a constant matrix, i.e.
\begin{equation}\label{mpc1}
\left(\frac{d\Psi}{d\xi}+ \frac{i}{2}\xi \sigma_3\Psi\right)\Psi^{-1}=\frac{i}{2}[\sigma_3,M_{1}^{pc}]=\left(\begin{array}{cc}
\beta_{11}^1 & \beta_{12}^1 \\
\beta_{21}^1 & \beta_{22}^1
\end{array}\right).
\end{equation}
Comparing both sides of the equation gives
\begin{equation}
(M_{1}^{pc})_{12}=-i\beta_{12}^1,~(M_{1}^{pc})_{21}=i\beta_{21}^1,~\beta_{11}^1=\beta_{22}^1=0.
\end{equation}
First consider the upper half-plane $Im\xi > 0$. Let\\
$~~~~~~~~~~~~~~~~~~~~~~~~~~~~~~~~~~~~~~~~~~~~~~\Psi=\left(\begin{array}{cc}
\Psi_{11} & \Psi_{12} \\
\Psi_{21} & \Psi_{22}
\end{array}\right),$\\
from \eqref{mpc1}, we get
\begin{equation}
\begin{aligned}
\partial_{\xi} \Psi_{11}^{+}+\frac{1}{2}i\xi \Psi_{11}^{+}=\beta_{12}^1 \Psi_{21}^{+}, \\
\partial_{\xi} \Psi_{21}^{+}-\frac{1}{2}i\xi \Psi_{21}^{+}=\beta_{21}^1 \Psi_{11}^{+},
\end{aligned}
\end{equation}
so have
\begin{equation}\label{w01}
\partial_{\xi}^2 \Psi_{11}^{+}=\left(-\frac{\xi^2}{4}-\frac{i}{2}+\beta_{12}^1\beta_{21}^1\right)\Psi_{11}^{+}.
\end{equation}
If we let $\Psi_{11}^{+}(\xi)=g(e^{-\frac{3\pi i}{4}}\xi)$, then \eqref{w01} reduces to the Weber equation
\begin{equation}
\partial_{\zeta}^2 g(\zeta)+\left(\frac{1}{2}-\frac{\zeta^2}{4}+a \right)g(\zeta)=0,
\end{equation}
where $a=i\beta_{12}^1\beta_{21}^1$. This is a second order ordinary differential equation with two linearly independent solutions $D_a(\zeta)$, $D_a(-\zeta)$. There are constants $c_1$, $c_2$ such that
\begin{equation}
\Psi_{11}^{+}=c_1D_a(e^{-\frac{3\pi i}{4}}\xi)+c_2D_a(-e^{-\frac{3\pi i}{4}}\xi),
\end{equation}
here $D_a(\zeta)$ is a standard parabolic function with the following asymptotics
\begin{equation}
D_{a}(\zeta)=\left\{\begin{array}{l}
\zeta^{a} e^{-\zeta^{2} / 4}\left(1+O\left(\zeta^{-2}\right)\right), \quad|\arg \zeta|<\frac{3 \pi}{4} \\
\zeta^{a} e^{-\zeta^{2} / 4}\left(1+O\left(\zeta^{-2}\right)\right)-(2 \pi)^{1 / 2} \Gamma^{-1}(-a) e^{a \pi i} \zeta^{-a-1} e^{\zeta^{2} / 4}\left(1+O\left(\zeta^{-2}\right)\right), \\
\frac{\pi}{4}<\arg \zeta<\frac{5 \pi}{4}, \\
\zeta^{a} e^{-\zeta^{2} / 4}\left(1+O\left(\zeta^{-2}\right)\right)-(2 \pi)^{1 / 2} \Gamma^{-1}(-a) e^{-a \pi i} \zeta^{-a-1} e^{\zeta^{2} / 4}\left(1+O\left(\zeta^{-2}\right)\right), \\
-\frac{5 \pi}{4}<\arg \zeta<-\frac{\pi}{4}.
\end{array}\right.
\end{equation}
By calculating, we get\\
$~~~~~~~~~~~~~~~~~~~~~~~~~~~~~~~~~~~~~~~\beta_{12}^1=\frac{(2\pi)^{\frac{1}{2}}e^{\frac{i\pi}{4}}e^{-\frac{\pi \nu}{2}}}{r_{z_1}\Gamma(-a)}$, $\beta_{12}^1\beta_{21}^1=\nu$,\\
$~~~~~~~~~~~~~~~~~~~~~~~~~~~~~~~~~~~~~~~\arg \beta_{12}^1=\frac{\pi}{4}-\arg r_{z_1}-\arg \Gamma(-i\nu)$.\\
By using the literature \cite{Deift93,Jenkins18}, we can obtain
\begin{equation}
M^{lo, 1}(z)=I+\frac{t^{-1 / 2}}{ \sqrt{2 \theta^{\prime \prime}\left(z_{1}\right)}\left(z-z_{1}\right)}\left(\begin{array}{cc}
0 & -i \beta_{12}^1 \\
i \beta_{21}^1 & 0
\end{array}\right)+\mathcal{O}\left(t^{-1}\right).
\end{equation}
In the same way, we get\\
$~~~~~~~~~~~~~~~~~~~~~~~~~~~~~~~~~~~~~~~\beta_{12}^2=\frac{(2\pi)^{\frac{1}{2}}e^{-\frac{i\pi}{4}}e^{\frac{\pi \nu}{2}}}{r_{z_2}\Gamma(a)}$, $\beta_{12}^2\beta_{21}^2=\nu$,\\
$~~~~~~~~~~~~~~~~~~~~~~~~~~~~~~~~~~~~~~~\arg \beta_{12}^1=-\frac{\pi}{4}-\arg r_{z_2}+\arg \Gamma(-i\nu)$,\\
where $r_{z_{2}}=r\left(z_{2}\right)e^{-2i\lambda(z_2,z_2)+2it\theta(z_2)}\exp(-i\nu(z_2)\log(2t\theta^{''}(z_2)))$.\\
And
\begin{equation}
M^{lo,2}(z)=I+\frac{t^{-1/2}}{ \sqrt{2 \theta^{''}\left(z_{2}\right)}\left(z-z_{2}\right)}\left(\begin{array}{cc}
0 & -i \beta_{12}^2 \\
i \beta_{21}^2 & 0
\end{array}\right)+\mathcal{O}\left(t^{-1}\right).
\end{equation}
Then combining Proposition 5.1, we obtain\\
$\mathbf{Proposition~ 5.2}$ As $t \rightarrow +\infty$,
\begin{equation}
M^{(lo)}(z)=I+\frac{t^{-1/2}}{\sqrt{2 \theta^{''}\left(z_{1}\right)}\left(z-z_{1}\right)}\left(\begin{array}{cc}
0 & -i \beta_{12}^1 \\
i \beta_{21}^1 & 0
\end{array}\right)
+\frac{t^{-1/2}}{\sqrt{2 \theta^{''}\left(z_{2}\right)}\left(z-z_{2}\right)}\left(\begin{array}{cc}
0 & -i \beta_{12}^2 \\
i \beta_{21}^2 & 0
\end{array}\right)+\mathcal{O}\left(t^{-1}\right).
\end{equation}

\section{Asymptotic analysis on $\bar{\partial}$-problem}

We use $M^{(lo)}(z)$ to construct a new matrix function
\begin{equation}\label{5h3}
M^{(3)}(z)=M^{(2)}(z)M^{(lo)}(z)^{-1},
\end{equation}
then $M^{(3)}(z)$ is continuous and has no jumps on $\mathbb{C}$, in fact,\\
$~~~~~~~~~~~~~~(M_{-}^{(3)}(z))^{-1}(z)M_{+}^{(3)}(z)=M_{-}^{(lo)}(z)(M_{-}^{(2)}(z))^{-1}M_{+}^{(2)}(z)(M_{+}^{(lo)}(z))^{-1}(z)$\\
$~~~~~~~~~~~~=M_{-}^{(lo)}(z)V^{(2)}(z)(M_{-}^{(lo)}(z)V^{(2)}(z))^{-1}=I$.\\
Thus we obtain a pure $\bar{\partial}$-problem:\\
\textbf{Pure $\bar{\partial}$-Problem:} Find $M^{(3)}(z)$ with following properties:\\
$\bullet $ Analyticity: $M^{3}(z)$ is continuous in $\mathbb{C}$;\\
$\bullet $ Asymptotic behavior:
\begin{equation}
M^{3}(z)\sim I+O(z^{-1}),~~z\rightarrow \infty;
\end{equation}
$\bullet $ $\bar{\partial}$-Derivative: We have\\
\begin{equation}
\bar{\partial}M^{3}(z)=M^{3}(z)W^{(3)}(z),~~z\in \mathbb{C};
\end{equation}
where
\begin{equation}
W^{(3)}(z)=M^{lo}(z)\bar{\partial}R^{(2)}(z)M^{lo}(z)^{-1}.
\end{equation}
The solutions of the $\bar{\partial}$-problem for $M^{3}(z)$ are equivalent to the integral equation
\begin{equation}\label{M(3)}
M^{(3)}(z)=I-\frac{1}{\pi} \iint_{\mathbb{C}} \frac{M^{(3)}(s) W^{(3)}(s)}{s-z} dA(s),
\end{equation}
where $A(s)$ is the Lebesgue measure on the $\mathbb{C}$. Then we can rewrite equation \eqref{M(3)} as
\begin{equation}\label{M3C}
(I-C_z)M^{(3)}(z)=I,
\end{equation}
where $C_z$ is the left Cauchy-Green integral operator
\begin{equation}
fC_{z}(z)=-\frac{1}{\pi} \iint_{\mathbb{C}} \frac{f(s) W^{(3)}(s)}{s-z} dA(s).
\end{equation}
$\mathbf{Proposition~ 6.1}$ For $t\rightarrow \infty$, the operator $C_z$ is a small norm, and
\begin{equation}
\| C_z\|_{L^{\infty}\rightarrow L^{\infty}} \lesssim t^{-1/4}.
\end{equation}
$\mathbf{Proof.}$ Taking the region $\Omega_{11}$ $(s=z_1+u+iv, z=x+iy)$ as an example, the
other regions are similarly to discuss. For any $f \in L_{\infty}$, we have
\begin{equation}
\begin{aligned}
\left\|f C_{z}\right\|_{L^{\infty}} &\leq \|f\|_{L^{\infty}} \frac{1}{\pi} \iint_{\Omega_{11}} \frac{\left|W^{(3)}(s)\right|}{|z-s|} dA(s),\\
&\leq \|f\|_{L^{\infty}} \frac{1}{\pi} \iint_{\Omega_{11}} \frac{\left|\bar{\partial}R^{(2)}(s)
\right|}{|z-s|} dA(s).
\end{aligned}
\end{equation}
Next, we consider only the integral equation
\begin{equation}
\iint_{\Omega_{11}} \frac{\left|\bar{\partial}R^{(2)}(s)\right|}{|z-s|} dA(s),
\end{equation}
From equation \eqref{z}, we know that when $z \in \Omega_{11}$
\begin{equation}
Im\theta(z)\leq -c|v||Rez-z_1|,
\end{equation}
we have
\begin{equation}
\begin{aligned}
&~~~\iint_{\Omega_{11}} \frac{\left|\bar{\partial}R^{(2)}(s)\right|}{|z-s|} dA(s)\\
&=\iint_{\Omega_{11}} \frac{|f_{11}^{'}(s)|e^{2tIm\theta}}{|s-z|}dA(s)+
\iint_{\Omega_{11}} \frac{|s-z_1|^{-\frac{1}{2}}e^{2tIm\theta}}{|s-z|}dA(s)\\
&=I_1+I_2,
\end{aligned}
\end{equation}
where $I_1=\iint_{\Omega_{11}} \frac{|f_{11}^{'}(s)|e^{2tIm\theta}}{|s-z|}dA(s)$, $I_2=\iint_{\Omega_{11}} \frac{|s-z_1|^{-\frac{1}{2}}e^{2tIm\theta}}{|s-z|}dA(s)$.

\begin{equation}
\begin{aligned}\label{s-z}
\left\|\frac{1}{s-z}\right\|_{L^{2}\left(v+z_{1}, \infty\right)}^{2} & =\int_{v+z_{1}}^{\infty} \frac{1}{|s-z|^{2}} d u \leqslant \int_{-\infty}^{\infty} \frac{1}{|s-z|^{2}} d u \\
& =\int_{-\infty}^{\infty} \frac{1}{(u-x)^{2}+(v-y)^{2}} d u \\
& =\frac{1}{|v-y|} \int_{-\infty}^{\infty} \frac{1}{1+\eta^{2}} d\eta=\frac{\pi}{|v-y|},
\end{aligned}
\end{equation}
where $\eta=\frac{u-x}{v-y}$.\\
Using \eqref{s-z}, the direct calculation gives
\begin{equation}
\begin{aligned}\label{I1}
I_{1} & \leq \int_{0}^{\infty} \int_{v}^{\infty} \frac{\left|f_{11}^{'}(s)\right| e^{-c t\left(Res-z_{1}\right) v}}{|s-z|} d u d v\\
& \leq \int_{0}^{\infty} e^{-ct v^{2}}\left\|f_{11}^{'}(s)\right\|_{L^{2}}\left\|\frac{1}{s-z}\right\|_{L^{2}} dv \\
 &\leq \int_{0}^{\infty} \frac{e^{-4t v^{2}}}{\sqrt{|v-y|}} dv\\
 & \leq t^{-1/4} .
\end{aligned}
\end{equation}

\begin{equation}\label{I2}
\begin{aligned}
I_{2} & \leq \int_{0}^{+\infty}e^{-ctv^2}dv\int_{v}^{+\infty}\frac{|s-z_1|^{-\frac{1}{2}}}{|s-z|}du\\
&\leq \int_{0}^{+\infty}e^{-ctv^2} \left\|\frac{1}{s-z}\right\|_{L^q(R^{+})} \left\||s-z_1|^{-\frac{1}{2}} \right\|_{L^p(R^{+})}dv\\
&=\left(\int_{0}^{y}+\int_{y}^{+\infty}\right)v^{-\frac{1}{2}+\frac{1}{p}}|v-y|^{\frac{1}{q}-1}e^{-ctv^2}dv.
\end{aligned}
\end{equation}
For the first integral, we have
\begin{equation}
\begin{aligned}
&~~~\int_{0}^{y}v^{-\frac{1}{2}+\frac{1}{p}}|v-y|^{\frac{1}{q}-1}e^{-ctv^2}dv\\
&=\int_{0}^{1}\sqrt{y}e^{-cty^2\omega^2}\omega^{\frac{1}{p}-\frac{1}{2}}|1-\omega|^{\frac{1}{q}-1} d\omega\\
&\leq t^{-\frac{1}{4}},
\end{aligned}
\end{equation}
for the second integral, we have
\begin{equation}
\begin{aligned}
&~~~\int_{y}^{+\infty}v^{-\frac{1}{2}+\frac{1}{p}}|v-y|^{\frac{1}{q}-1}e^{-ctv^2}dv\\
&=\int_{0}^{+\infty}e^{-t(y+\omega)^2}(y+\omega)^{\frac{1}{p}-\frac{1}{2}}\omega^{\frac{1}{q}-1} d\omega\\
&\leq \int_{0}^{+\infty}e^{-ct\omega^2}\omega^{\frac{1}{p}-\frac{1}{2}}\omega^{\frac{1}{q}-1}d\omega\\
&\leq \int_{0}^{+\infty}e^{-t\omega^2}\omega^{-\frac{1}{2}}d\omega\\
&\leq t^{-\frac{1}{4}}.
\end{aligned}
\end{equation}
For $z \in \Omega_{jk},j=1,2;k=1,2,3,4$, we can conclude that $\|C_z\|_{L^{\infty}(\mathbb{C})\rightarrow L^{\infty}(\mathbb{C})} \leq t^{-\frac{1}{4}}$.

\section{The long-time asymptotics of the NHNSE}
$M^{(3)}(z)$ is expanded as follows
\begin{equation}
M^{(3)}(z)=I+\frac{M_{1}^{(3)}(z)}{z}+
\frac{1}{\pi} \iint_{\mathbb{C}} \frac{sM^{(3)}(s) W^{(3)}(s)}{z(z-s)} dA(s),
\end{equation}
where\\
$~~~~~~~~~~~~~~~~~~~~~~~~~~~~M_{1}^{(3)}(z)=\frac{1}{\pi} \iint_{\mathbb{C}} M^{(3)}(s) W^{(3)}(s)dA(s)$.\\
$\mathbf{Proposition~ 7.1}$ There exists a constant $c$, such that
\begin{equation}
|M_{1}^{(3)}(z)|\leq ct^{-\frac{3}{4}}.
\end{equation}
$\mathbf{Proof}$
\begin{equation}
\begin{aligned}
|M_{1}^{(3)}(z)|& \leq \frac{1}{\pi} \iint_{\Omega_{11}}\left|M^{(3)}(s)M^{lo}(s)\bar{\partial}R^{(2)}(s)M^{lo}(s)^{-1} \right|dA(s)\\
&\leq \frac{1}{\pi} \iint_{\Omega_{11}}\left|M^{(3)}(s)\right| \left|M^{lo}(s)\right| \left|\bar{\partial}R^{(2)}(s)\right| \left|M^{lo}(s)^{-1} \right|dA(s)\\
&\leq c(I_3+I_4),
\end{aligned}
\end{equation}
where\\
$~~~~~~~~~~~~~~~~~~~~~~~~~~~~~~~~~~~~~~~~I_3=\iint_{\Omega_{11}}|f_{11}^{'}(s)|e^{2tIm\theta}dA(s)$,\\
$~~~~~~~~~~~~~~~~~~~~~~~~~~~~~~~~~~~~~~~~I_4=\iint_{\Omega_{11}}|s-z_1|^{-\frac{1}{2}}e^{2tIm\theta}dA(s)$.\\
Using the Cauchy-Schwarz inequality,
\begin{equation}
\begin{aligned}
|I_3|&\leq c \iint_{\Omega_{11}}e^{-2tv(Res-z_1)}dA(s)\\
&\leq c \int_{0}^{+\infty}\left(\int_{v}^{+\infty} e^{-4tuv}du \right)^{\frac{1}{2}}dv\\
&\leq c t^{-\frac{1}{2}}\int_{0}^{+\infty}v^{-\frac{1}{2}}e^{-2tv^2}dv\\
&\leq ct^{-\frac{3}{4}}.
\end{aligned}
\end{equation}
Via Holder inequality,
\begin{equation}
\begin{aligned}
|I_4|& \leq c \int_{0}^{+\infty}v^{\frac{1}{p}-\frac{1}{2}} \left(\int_{v}^{+\infty} e^{-4tquv}du \right)^{\frac{1}{q}}dv\\
& \leq c \int_{0}^{+\infty}v^{\frac{1}{p}-\frac{1}{2}} (qtv)^{-\frac{1}{q}} e^{-2tv^2}dv\\
&\leq t^{-\frac{1}{q}}\int_{0}^{+\infty}v^{\frac{2}{p}-\frac{3}{2}}e^{-2tv^2}dv\\
&\leq  t^{-\frac{3}{4}}\int_{0}^{+\infty}\omega^{\frac{2}{p}-\frac{3}{2}}e^{-2\omega^2}d\omega\\
&\leq ct^{-\frac{3}{4}}.
\end{aligned}
\end{equation}

Review the series of transformations \eqref{t1}, \eqref{t2} and \eqref{t3}, we have
\begin{equation}
M(z)=M^{(3)}(z)M^{(lo)}(z)R^{(2)}(z)^{-1}\delta(z)^{\sigma_3}.
\end{equation}
To construct the solution $q(x, t)$, we take $z\rightarrow \infty$ along the imaginary axis. The
advantage of this is that $R^{(2)}(z)$ is the identity matrix, so
\begin{equation}
M(z)=\left(I+\frac{M_1^{(3)}(z)}{z}+\cdots\right) \left(I+\frac{M_1^{(lo)}(z)}{z}+\cdots \right)\delta(z)^{\sigma_3},
\end{equation}
from which, comparing the coefficients of $z^{-1}$ to get
\begin{equation}
M_1(z)=I+M_1^{(3)}(z)+M_1^{(lo)}(z),
\end{equation}
where\\
$~~~~~~~~~~~~~~~~~~M_1^{(lo)}(z)=\frac{t^{-1/2}}{\sqrt{2 \theta^{''}\left(z_{1}\right)}}\left(\begin{array}{cc}
0 & -i \beta_{12}^1 \\
i \beta_{21}^1 & 0
\end{array}\right)
+\frac{t^{-1/2}}{\sqrt{2 \theta^{''}\left(z_{2}\right)}}\left(\begin{array}{cc}
0 & -i \beta_{12}^2 \\
i \beta_{21}^2 & 0
\end{array}\right)$.\\
Then from \eqref{q}, we get
\begin{equation}
q(x,t)=2i(M_1^{(lo)}(z))_{12}+O(t^{-\frac{3}{4}}),
\end{equation}
which is a novel result.

\section{Conclusion}
A new Lax pair was introduced whose compatibility condition lead to the well-known KdV equation, the modified KdV equation and the standard nonlinear Schr\"{o}dinger equation. In particular, a new high-order nonlinear Schr\"{o}dinger equation (NHNSE) was derived from the Lax pair. In the paper, we focused on the study of some properties of the NHNSE, including the analysis of the characteristic function, the constructing the RH problem, the continuous extension of the scattering data, the long-time asymptotic behavior. It is remarkable that such the properties mentioned as above have not been investigated by other research people. We obtained completely novel long-time asymptotic behavior of the NHNSE.

\section*{Acknowledgements}

This work was supported by the National Natural Science Foundation of China grant No.12371256;
the National Natural Science Foundation of China grant No.11971475; SuQian Sci\&Tech Program Grant No.K202225.\\


\end{CJK*}
\end{document}